\documentclass[12pt]{amsart}
\usepackage{amsfonts,amsmath,amssymb,amscd,amsthm,enumerate}
\usepackage{amsrefs}
\usepackage{geometry}
\newtheorem{theorem}{Theorem}
\newtheorem{prop}[theorem]{Proposition}
\newtheorem{lem}[theorem]{Lemma}
\newtheorem{cor}[theorem]{Corollary}
\theoremstyle{definition}

\newtheorem{rem}[theorem]{Remark}
\newtheorem{mydef}[theorem]{Definition}
\newtheorem{example}[theorem]{Example}

\newtheorem{notation}[theorem]{Notation}

\def\C{\mathbb{C}}
\def\L{\mathbb{L}}
\def\R{\mathbb{R}}

\def\T{\mathbb{T}}
\def\E{\mathbb{E}}
\def\I{\mathbb{I}}

\def\x{\mathbf{x}}

\def\v{\mathbf{v}}
\def\u{\mathbf{u}}
\def\y{\mathbf{y}}
\def\z{\mathbf{z}}
\def\w{\mathbf{w}}
\def\e{\mathbf{e}}
\def\<{\langle}
\def\>{\rangle}
\def\Span{{\rm span}}
\def\SGF{{\rm SGF}}
\def\SGF{{\rm SGF}}
\def\CRF{{\rm CRF}}

\begin{document}
\title{An Abstract Morimoto Theorem for Generalized $F$-structures}
\author{Marco Aldi and Daniele Grandini}

\begin{abstract}
We abstract Morimoto's construction of complex structures on product manifolds to pairs of certain generalized $F$-structures on manifolds that are not necessarily global products. As an application we characterize invariant generalized complex structures on products in which one factor is a Lie group and generalize a theorem of Blair, Ludden and Yano on Hermitian bicontact manifolds.
\end{abstract}
\maketitle
\section{Introduction}
The study of generalized geometry in arbitrary (not necessarily even) dimension was pioneered by Vaisman \cite{V} and further developed by various authors (\cite{PW},\cite{S},\cite{GT1},\cite{AG},\cite{GT2}). The key notion is that of {\it generalized $F$-structure} i.e.\ a skew-symmetric endomorphism $\Phi:\T M\rightarrow \T M$ of the generalized tangent bundle $\T M=TM\oplus T^*M$ of a manifold M, such that $\Phi^3+\Phi=0$. It easy to see that if $\Phi$ is a generalized $F$-structure on $M$, then the restriction of the tautological inner product to the kernel of $\Phi$ is nondegenerate on each fiber. In this paper we focus on a specific kind of generalized $F$-structures, for which $\ker(\Phi)$ has fiberwise split signature. Most natural examples of generalized $F$-structures, including generalized almost complex structures and generalized almost contact structures, have split signature. To study generalized $F$-structures, we find it convenient to first introduce the notion of \emph{split structure} i.e.\ a subbundle $E\subseteq \T M$ on which the tautological inner product is nondegenerate and has split signature. A {\it split generalized $F$-structure} (or SGF-structure) is then defined to be an orthogonal, skew-symmetric endomorphism $J$ of a split structure $E$.

The generalized tangent bundle is acted upon by the group ${\rm Diff}(M)\ltimes \Omega_{cl}(M)$ of extended diffeomorphism with closed forms acting by the so-called B-field transform. Infinitesimally, this action corresponds to the notion of {\it generalized Lie derivative} $\L_{\x}$ \cite{H}. Given a subset $S\subseteq \Gamma(\T M\otimes \C)$, it is useful to consider its {\it normalizer} $\I(S)$ i.e.\ the set of all sections $\x$ of $\T M\otimes \C$ such that $\L_{\x}(S)\subseteq S$. By definition, the normalizer of a split generalized $F$-structure $J\in {\rm End}(E)$ is the normalizer $\I(J)$ of its $\sqrt{-1}$-eigenbundle $L_J$. Geometrically, $\I(J)$ can be thought of the set of infinitesimal symmetries of $E$ that commute with $J$. According to \cite{V}, $J$ is a {\it generalized \CRF-structure} if $L_J\subseteq \I(J)$.

Given two SGF-structures $J_1$, $J_2$ it is natural to ask under what conditions $L_{J_1}$ normalizes $L_{J_2}$. For instance, if $(M_1,\mathcal J_1)$ and $(M_2,\mathcal J_2)$ are generalized almost complex structures, then $\mathcal J_1$ (resp.\ $\mathcal J_2$) lifts to an SGF-structure $J_1$ (resp.\ $J_2$) on the split structure  on $M_1\times M_2$ generated by sections of $\T M_1$ (resp.\ of $\T M_2$). It is then easy to see that $L_{J_1}$ and $L_{J_2}$ normalize each other and that $J_1\oplus J_2$ is a generalized complex structure if and only if both $J_1$ and $J_2$ are. Similarly, if $\mathcal J_1$ and $\mathcal J_2$ are SGF-structures corresponding to generalized almost contact structures on $M_1$ and $M_2$, one can still define their lifts $J_1$, $J_2$ to $M_1\times M_2$, but $J_1\oplus J_2$ is no longer a generalized almost complex structure for dimensional reasons. However, generalizing a classical construction of Morimoto, one can introduce a third SGF-structure $\Psi$ on $M_1\times M_2$ in such a way that $J_1\oplus J_2\oplus \Psi$ is a generalized almost complex structure. Extending a theorem of Morimoto \cite{M} to the generalized setting, Gomez and Talvacchia \cite{GT1} proved the existence of a canonical \SGF-structure $\Psi$ for which $J_1\oplus J_2\oplus \Psi$ is a generalized complex structure if and only if $J_1$, $J_2$ and $\Psi$ are generalized \CRF-structures and the natural framing of $L_\Psi\oplus \overline L_\Psi$ normalizes both $J_1$ and $J_2$.

In this paper we abstract the features that make Morimoto's construction \cite{M} work into the concept of {\it (adaptable) Morimoto datum} defined out of: 1) mutually orthogonal split structures $E_1$, $E_2$, $E_1'$ and $E_2'$; 2) SGF-structures $J_1$ on $E_1$, $J_2$ on $E_2$, $\Psi$ on $E_1'\oplus E_2'$ and 3) global framings $V_1$ for $E_1'$ and $V_2$ for $E_2'$. Our main results is an {\it Abstract Morimoto Theorem} stating that in presence of an adaptable Morimoto datum, $J_1\oplus J_2\oplus \Psi$ is a generalized \CRF-structure if and only if $(J_1,V_1)$ and $(J_2,V_2)$ are {\it normal pairs}, an abstraction of the concept of normal generalized almost contact structure introduced in \cite{AG}.

Our Abstract Morimoto Theorem unifies and extends several theorems \`a la Morimoto in the literature. If $M$ is indeed a product $M_1\times M_2$ and $E_i$, $E_i'$ are pull-back of split structures on $\T M_i$, then our construction yields generalized almost complex structures on $M_1\times M_2$ which simultaneously generalize Morimoto products of generalized almost contact structures \cite{GT1} and Morimoto products of classical framed $F$-structures \cite{N}. The generalized complex structures constructed with our method come in families and thus, even in the generalized contact case, they are more general than those of \cite{GT1}. For instance, we show that the Morimoto product of two copies of the normal generalized almost contact structures on $S^3$ introduced in \cite{AG} yields holomorphic Poisson deformations of the Calabi-Eckmann complex structures on $S^3\times S^3$ for every choice of complex structure on the $T^2$ fiber. In a different direction, we are able to extend Sekiya's characterization of invariant generalized (almost) complex structures (\cite{S}, \cite{AG}) from products of the form $M\times \R$ to products of $M$ with an arbitrary finite dimensional Lie group.

An important feature of the notion of Morimoto datum is that it is sufficiently flexible to apply to manifolds that are not necessarily global products. For instance, we are able to describe two constructions of generalized \CRF-structures on flat principal bundles, one of which extends previous work \cite{BH} on normal contact pairs. A second class of examples of Morimoto data beyond the global product case  comes from a generalized version of a classical theorem of Blair, Ludden and Yano \cite{BLY} which states that Hermitian bicontact manifolds with bicontact forms $(\eta_1,\eta_2)$ of bidegree $(1,1)$ are locally the product of normal contact manifolds. In this paper we prove an Abstract Blair-Ludden-Yano Theorem at the level of Hermitian bicontact data, a notion that we introduce in order to isolate the features of classical Hermitian bicontact structures of bidegree $(1,1)$ that we need. On the one hand, we prove that our Abstract Blair-Ludden-Yano Theorem implies the classical one. On the other hand, we show that this generalization is non-trivial since the non-commutative Calabi-Eckmann structures on $S^3\times S^3$ provide non-classical examples of Hermitian bicontact data.

The paper is organized as follows. Section \ref{sec:preliminaries} is a recollection of basic notions and notations used in generalized geometry. We refer the reader to \cite{G} and \cite{H} for a systematic treatment of the subject. In Section \ref{sec:split} we define our main objects of study: split structures, SGF-structures and split generalized \CRF-structures. In Section \ref{sec:normal} we study normalizers of SGF-structures and introduce the important notion of normal pair. Section \ref{sec:morimoto} contains the definition of Morimoto datum and the Abstract Morimoto Theorem. Section \ref{sec:connection} is technical in nature and describes the behavior of normalizers and normal pairs under pull-back by a surjective submersion. In Section \ref{sec:externalmorimoto}, Section \ref{sec:sekiya} and Section \ref{sec:flat} we specialize the Abstract Morimoto Theorem to various particular cases including global products and flat principal bundles, making the connection with previous results in the literature. We conclude with Section \ref{sec:BLY} in which we introduce the concept of Hermitian bicontact datum and prove the Abstract Blair-Ludden-Yano Theorem. In this paper, the notion of contact and bicontact datum is developed mainly for the purpose of providing non-trivial examples of Morimoto data. A systematic treatment of (bi)contact data, in particular exploring their connection with other attempts to extends contact geometry to the generalized setting (e.g.\ \cite{PW}, \cite{IW}), would be interesting and we hope to come back to this point in the future.

\section{Preliminaries on Generalized Geometry}\label{sec:preliminaries}

\begin{mydef}
The \emph{generalized tangent bundle} of a real smooth manifold $M$ of finite dimension $n$ is the vector bundle $\T M:=TM\oplus T^*M$. $\T M$ is endowed with a $C^{\infty}(M)$-bilinear, symmetric {\it tautological inner product} of signature $(n,n)$ defined by
\[
\langle X+\alpha, Y+\beta\rangle:=\frac{1}{2}(\alpha(Y)+\beta(X))
\]
for all $X,Y\in \Gamma(TM)$ and all $\alpha, \beta\in \Gamma(T^*M)$. The generalized tangent bundle is also endowed with an $\R$-bilinear map $[\ ,\ ]: \Gamma(\T M)\times \Gamma(\T M)\rightarrow \Gamma(\T M)$ called the \emph{Dorfman bracket}, defined by
\[
[X+\alpha, Y+\beta]:=[X,Y]+{\mathcal L}_X\beta-\iota_Yd\alpha
\]
for all $X,Y\in \Gamma(TM)$ and for all $\alpha, \beta\in \Gamma(T^*M)$. Sections of $\T M$ are denoted by $\x, \y$, etc. unless their (co)tangent components need to be specified.
\end{mydef}

\begin{mydef}
For each $\x \in \Gamma(\T M)$, the {\it generalized Lie derivative with respect to $\x$} is the $\R$-linear endomorphism $\L_\x$ of $\Gamma(\T M)$ defined by $\L_\x(\y)=[\x,\y]$ for all $\y\in \Gamma(\T M)$. $\L_\x$ extends to the unique endomorphism of the full tensor algebra of $\Gamma(\T M)$ such that $\L_\x(f)=2\langle \x,df\rangle$ for all $f\in C^\infty(M)$ and such that $\L_x$ is a graded derivation with respect to the tensor product.
\end{mydef}

\begin{rem}
Let $a$ be the projection of $\T M$ onto the tangent bundle $TM$. The quadruple $(\T M, \langle\ ,\ \rangle, [\ ,\ ],a)$ satisfies the axioms of \emph{Courant algebroid}
\begin{enumerate}[i)]
  \item $a(\x)\left(\langle\y,\z\rangle\right)=\langle[\x,\y],\z\rangle+\langle\y,[\x,\z]\rangle$,
  \item $[\x,[\y,\z]]=[[\x,\y],\z]+[\y,[\x,\z]]$,
  \item $[\x,\y]+[\y,\x]=2d\langle\x,\y\rangle$,
\end{enumerate}
for all $\x,\y,\z\in \Gamma(\T M)$. These  properties can be restated in terms of generalized Lie derivatives as follows
\begin{align}
\L_\x \langle \y,\z \rangle &= \langle \L_\x (\y),\x\rangle +\langle \y, \L_\x(\z)\rangle\,;\\
\L_\x[\y,\z] &= [\L_\x(\y),\z]+[\y,\L_\x(\z)]\,;\\
2d\langle \x,\y\rangle&=\L_\x(\y)+ \L_\y(\x)\,;
\end{align}
for all $\x,\y,\z\in \Gamma(\T M)$.
\end{rem}

\begin{rem} It is well-known that given a closed three-form $H$ on $M$, one may twist the Dorfman bracket to
\[
[\x,\y]_{H}:=[\x,\y]-\iota_{\x}\iota_{\y}H
\]
which also satisfies the axioms of Courant algebroid. While the results of this paper rely only on these and therefore extend to the twisted case, we set $H=0$ for notational convenience. 
\end{rem}

\begin{notation}Given a subset $S$ of $\Gamma(\T M)$, we denote the $C^{\infty}(M)$-submodule of $\Gamma(\T M)$ generated by $S$ by $\Span(S)$. We reserve the notation $\Span_{\R}(S)$ for the $\R$-submodule of $\Gamma(M)$ generated by $S$.
\end{notation}

\begin{mydef} 
Let $E$ be a subbundle of $\Gamma(\T M)$.  A \emph{framing} of $E$ is a real subspace $V$ of $\Gamma(E)$ whose dimension equals the rank of $E$ and such that $\Span(V)=\Gamma(E)$. Moreover, if $U$ is an open set in $M$, a \emph{local framing} of $E$ on $U$ is a framing of $E_{|U}$.
\end{mydef}

\section{Split structures}\label{sec:split}
\begin{mydef}Let $M$ be an $n$-dimensional manifold. A \emph{split structure} on $M$ of rank $2k$ is a subbundle $E\subseteq \T M$ such that the restriction $\langle,\rangle_{|E}$ is nondegenerate with signature $(k,k)$. We denote by $\E_k(M)$ the set of all split structures of rank $2k$ on $M$, and we write $\E(M)$ for the set of all split structures on $M$.
\end{mydef}

\begin{rem}\label{rem:split}
Split structures are closed with respect to the following operations.
\begin{enumerate}
\item If $E\in \E_k(M)$, then
\[
E^{\bot}=\{\x\in \Gamma(\T M)\,|\,\langle \x, \y\rangle=0\,\text{for all}\,\, \y\in E\}
\]
is a split structure of rank $2n-2k$.
\item Let $E\in \E(M)$, let $F:E\rightarrow \T M$ be a base preserving morphism and let $C$ be a nowhere vanishing function on $M$ such that
\[
\langle F\x,F\y \rangle=C\langle \x,\y\rangle
\]
for all $\x,\y\in \Gamma(E)$. Then $F(E)\in \E(M)$.
\item If $E\in \E_k(M)$ and $E'\in \E_{k'}(M)$ are such that $\langle \Gamma(E),\Gamma(E')\rangle=0$, then the Whitney sum $E \oplus E'$ is in $\E_{k+k'}(M)$.
\item If $E\in \E_k(M)$ and $E'\in \E_{k'}(M')$ then the external Whitney sum $E \boxplus E'$ is in $\E_{k+k'}(M\times M')$. Note that the space of sections $\Gamma(E)$ (resp.\ $\Gamma(E')$) is included canonically into the space $\Gamma(E \boxplus E')$ as a $C^{\infty}(M)$-submodule (resp.\ $C^{\infty}(M')$-submodule), but not a $C^{\infty}(M\times M')$-submodule.
\end{enumerate}
\end{rem}

\begin{rem}
If $E\in \E_k(M)$ is equipped with a framing $V$,then  the restriction of the tautological inner product to $V$ is nondegenerate with signature $(k,k)$. Moreover, the orthogonal group ${\rm O}(V)\subseteq O(E)$ can be identified (as a Lie group) with the subgroup endomorphisms $\Psi$ such that $\Psi(V)\subseteq V$.
\end{rem}

\begin{mydef}
Let $E\in \E(M)$. A \emph{split generalized $F$-structure on $E$} is a bundle endomorphism $J\in {\rm End}(E)$ which is skew-symmetric and orthogonal with respect to the tautological inner product. We denote by ${\rm SGF}(E)$ the set of all almost complex split structures on $E$.
\end{mydef}

\begin{rem}\label{rem:Vaisman}
Split generalized $F$-structures are a particular case of the generalized $F$-structures introduced in \cite{V}. In particular, the following two characterizations of $\SGF(E)$ can be easily deduced from the results of \cite{V}. Extending $J\in \SGF(E)$ to $\T M$ by $0$ provides a bijection between $\SGF(E)$ and the set of all orthogonal endomorphisms $\Phi$ of $\T M$ such that $\Phi^3+\Phi=0$ and $\ker(\Phi)=E$. On the other hand, assigning to $J$ the subbundle
\[
L_J=\{\x-\sqrt{-1}J{\x}\,|\,\x\in E\}
\]
defines a bijection between $\SGF(E)$ and the set of maximally isotropic subbundles $L$ of $E\otimes \C$ such that $L\cap \overline{L}=0$.
\end{rem}

\begin{example}
Viewing $\T M$ as split structure on $M$, $\SGF(\T M)$ coincides with the set of all generalized almost complex structures on $M$, as defined in \cite{G}.
\end{example}

\begin{example}\label{ex:GACS}
In \cite{AG}, a {\it generalized almost contact structure} is defined as a pair $(E,L)$ where $E\in \E_1(M)$ is a trivial subbundle of $\T M$ and $L$ is a maximal isotropic subbundle of $E^\perp\otimes \C$ such that $L\cap \overline L=0$. By Remark \ref{rem:Vaisman}, for each trivial $E\in \E_1(M)$ there is a canonical bijection between $\SGF(E^\perp)$ and the set of generalized almost contact structures of the form $(E,L)$. Let $J$ be the split generalized $F$-structure on $E^\perp$ corresponding to a generalized almost contact structure $(E,L)$ and let $\Phi$ be the extension of $J$ to $\T M$ by $0$. Given an isotropic frame $\{\e_1,\e_2\}$ of $E$ such that $2\langle \e_1,\e_2\rangle=1$ then $(\Phi,\e_1,\e_2)$ is a \emph{generalized almost contact triple} as defined in \cite{AG}. Therefore, the set of generalized contact triples up to a change of frame of $E$ can be identified with the union of all $\SGF(E^{\bot})$, as $E$ ranges over all rank 2 split structures on $M$ that are trivial subbundles of $\T M$.
\end{example}

\begin{example}\label{ex:fstructures}
If $\Phi$ is a classical $F$-structure in the sense of \cite{V}, then $\ker(\Phi)\cap TM$ and $\ker(\Phi)\cap T^*M$ are maximally isotropic in $\ker (\Phi)$. Therefore, the restriction of $\Phi$ to the orthogonal complement of $\ker(\Phi)$ is a split generalized $F$-structure.
\end{example}

\begin{mydef}
A split generalized $F$-structure $J\in \SGF(E)$ is a {\it split generalized \CRF-structure on $E\in \E(M)$} if its $\sqrt{-1}$-eigenbundle $L_J$ is closed under the Dorfman bracket. We denote by $\CRF(E)$ the set of all split generalized \CRF-structures on $E$.
\end{mydef}

\begin{example}
The set of all generalized complex structures on $M$ coincides with $\CRF(\T M)$.
\end{example}

\begin{example}\label{ex:S^3}
The following family of generalized almost contact structures on $M=S^3$ found in \cite{AG} will serve as a recurring example to illustrate the scope of the methods introduced in the present paper. Let $\{X_1,X_2,X_3\}$ be a global frame of $TS^3$ with dual frame $\{\alpha_1,\alpha_2,\alpha_3\}\subseteq T^*S^3$ such that $[X_i,X_j]= 2\varepsilon_{ijk} X_k$ and $[X_i,\alpha_j]=2\varepsilon_{ijk} \alpha_k$, where $\varepsilon_{ijk}$ is the Levi-Civita symbol. Given $h=f_2+\sqrt{-1}f_3\in C^\infty(S^3,\C)$, we deform $\alpha_1,\alpha_2,\alpha_3$ in the generalized sense to
\begin{align*}
\x_1&=\alpha_1+f_2X_2+f_3X_3\,,\\
\x_2&=\alpha_2-f_2X_1\,,\\
\x_3&=\alpha_3-f_3X_1\,.
\end{align*}
This leads to an interesting decomposition of $\T S^3$ as orthogonal direct sum of the split structures  $E=\Span(X_2,X_3,\x_2,\x_3)$ and $E'=\Span(X_1,\x_1)$. For any $h$, we also consider the split generalized $F$-structure $J\in\SGF(E)$ defined by $J(X_2)=X_3$ and $J(\x_2)=\x_3$. If $h=0$, we recover the standard almost contact structure on $S^3$ written in coordinates for which $X_1$ is tangent to the fibers of the Hopf fibration. A direct calculation shows that $J\in \CRF(E)$ if and only if $\overline\partial(h)=0$, where $\partial=X_2-\sqrt{-1}X_3$.
\end{example}

\section{Normalizers and Normal Pairs}\label{sec:normal}

\begin{mydef} Let $S$ be a subset of $\Gamma(\T M \otimes \C)$. We say that a section $\x$ of $\T M\otimes \C$ {\it normalizes $S$} if $\L_\x(S)\subseteq S$. The set $\I(S)$ of all sections that normalize $S$ is called the {\it normalizer} of $S$. If $T\subseteq \T M\otimes \C$ is a subbundle, we simply write $\I(T)$ for $\I(\Gamma(T))$.
\end{mydef}

\begin{rem}\label{rem:normalizer}
Let $E\in \E(M)$. Given $\x\in \I(E)$, $\y\in
\Gamma(E^\perp)$ and $\z\in \Gamma(E)$,
\[
0=\L_\x\langle \y,\z\rangle = \langle \L_\x(\y),\z \rangle + \langle \y,\L_\x(\z)\rangle = \langle \L_\x(\y),\z \rangle
\]
from which we conclude that $\I(E)=\I(E^\perp)$.
\end{rem}

\begin{mydef}
If $J$ is a split generalized $F$-structure and $L_J$ is its $\sqrt{-1}$-eigenbundle, we define the {\it normalizer of $J$} to be $\I(J)=\I(L_J)$. Given two split generalized $F$ structures $J_1$ and $J_2$ we say that $J_1$ {\it normalizes} $J_2$  if $\Gamma(L_{J_1})\subseteq \I(J_2)$.
\end{mydef}

\begin{example}
Let $J$ be a split generalized $F$-structure on $E\in \E(M)$. Then $J\in \CRF(E)$ if and only if $\Gamma(L_J)\subseteq \I(J)$.
\end{example}

\begin{rem} Let $E\in \E(M)$ and $J\in \SGF(E)$. Then $\x\in\I(J)$ if and only if $\x\in\I(E)$ and $\L_\x$ commutes with $J$ as elements of ${\rm End}_{\R}(\Gamma(E))$. By extending the action of $\L_\x$ to ${\rm End}_{\R}(\Gamma(E))$, this last requirement can we rewritten as $\L_{\x}(J)=0$.
\end{rem}

\begin{example}
Consider a generalized almost contact triple $(\Phi,\e_1,\e_2)$ as in Example \ref{ex:GACS}, let $E=\ker (\Phi)$ and let $J$ be the restriction of $\Phi$ to $E^\perp$. In the language of \cite{AG}, if $(\Phi,\e_1,\e_2)$ is integrable (resp.\ strongly integrable) then $J\in\SGF(E)$ is normalized by at least one of (resp.\ both) $\e_1$ and $\e_2$.
\end{example}

\begin{lem}\label{lem:Ju}
Let $J$ be a split generalized \CRF-structure on $E\in \E(M)$ and let $\u\in \I(J)$. Then $J(\u)\in \I(J)$.
\end{lem}

\noindent\emph{Proof:} Let $\v=J(\u)$. For every $\x\in \Gamma(E)$,
\[
[\u-\sqrt{-1}\v,\x-\sqrt{-1}J(\x)]=[\u,\x]-[\v,J(\x)]-\sqrt{-1}([\u,J(\x)]+[\v,\x])\,.
\]
Since $J\in \CRF(E)$ and $\u\in \I(J)$, then
\[
[\u,J(\x)]+[\v,\x]=J([\u,\x]-[\v,J(\x)])=[\u,J(\x)]-J[\v,J(\x)]\, ,
\]
which in turn implies $\v\in \I(J)$.
\qed

\begin{rem}
Due to the local nature of the Dorfman bracket, the normalizer of a subbundle $S\subseteq\T M$ defines a sheaf on $M$, whose sections on an open set $U\subseteq M$ are given by
\[
\I_U(S):=\{\x\in \Gamma_U(\T M):\L_{\x}\Gamma_U(S)\subseteq \Gamma_U(S)\}\,.
\]
\end{rem}

\begin{mydef}\label{def:complete}
A split structure $E\in \E(M)$ is said to be \emph{complete} if $\Gamma(E)$ is locally generated by $\Gamma(E)\cap\I(E)$, i.e.\ if each $p\in M$ admits an open neighborhood $U$ and a local framing  $W_U$ of $E$ on $U$, such that  $W_U\subseteq \Gamma_U(E)\cap \I_U(E)$.
\end{mydef}

\begin{mydef}
Let $E,E'\in \E(M)$ be such that $E'\subseteq E^\perp$, let $J\in\CRF(E)$ and let $V$ be a framing of $E'$. We say that $(J,V)$ is a {\it normal pair} if $V\subseteq \I(J)\cap \I(E')$.
\end{mydef}

\begin{example}
If $J\in\SGF(\T M)$, then $(J,0)$ is a normal pair if and only if $J$ is a generalized complex structure.
\end{example}

\begin{example}\label{ex:S^3normal}
Let $E,E'$ and $J$ be as in Example \ref{ex:S^3} and consider the framing $V=\Span_\R(X_1,\x_1)$ of $E'$. Then $(J,V)$ is a normal pair if and only if $h$ is annihilated by both $\overline\partial$ and $Y_1=X_1+2\sqrt{-1}{\rm Id}$.
\end{example}

\begin{example}\label{ex:normalGACS}
More generally, let $(\Phi,\e_1,\e_2)$ be a generalized contact triple as in Example \ref{ex:GACS}. Consider the framing $V={\rm span}_\R(\e_1,\e_2)$ of $E=\ker \Phi$ and denote by $J$ the restriction of $\Phi$ to $E$. Then $(J,V)$ is a normal pair if and only if $(\Phi,\e_1,\e_2)$ is a normal generalized contact triple in the sense of \cite{AG}. In this case, the condition $V\subseteq\I(E)$ implies that the Dorfman bracket vanishes identically on $V$.
\end{example}

\begin{lem}\label{lem:normality}
Let $E,E'\in \E(M)$ be such that $E'\subseteq E^\perp$. Given $J\in\CRF(E)$ and a framing $V$ of $E'$, the following are equivalent:
\begin{enumerate}[i)]
\item $(J,V)$ is a normal pair;
\item $V\subseteq \I(J)\cap \I(E\oplus E')$;
\item $V\subseteq \I(E) \cap \I(E\oplus E')$;
\item $V\subseteq \I(E) \cap \I(E')$.
\end{enumerate}
In particular, if $E'=E^\perp$, then $(J,V)$ is a normal pair if and only if $V\subseteq \I(E)$.
\end{lem}

\noindent\emph{Proof:} $V\subseteq \I(J)$ is equivalent to $\L_\v(\Gamma(L))\subseteq \Gamma(L)$ for all $\v\in V$, where $L$ is the $\sqrt{-1}$-eigenspace of $J$. Since $V=\overline V$ and $E\otimes \C=L\oplus \overline L$, this implies $V\subseteq \I(E)$. Therefore, i)$\Rightarrow$ii)$\Rightarrow$iii). If iii) holds, then for every $\v\in V$, $\e\in \Gamma(E)$ and $\e'\in \Gamma(E^{\bot})$

\[
\langle \L_\v\e',\e\rangle = \L_\v\langle \e',\e\rangle - \langle \e', \L_\v\e\rangle = 0
\]
which in turns implies iv). Under the assumptions of iv), $\L_\v(\Gamma(L))\subseteq \Gamma(E)$ for each $\v\in V$ and thus
\[
\langle \L_\v\x,\y \rangle = -\langle \L_\x\v,\y \rangle=-\L_\x\langle \v,\y\rangle+\langle \v,\L_\x\y\rangle =0
\]
for every $\x,\y\in \Gamma(L)$. Therefore, $\L_\v(\Gamma(L))\subseteq \Gamma(E\cap L^\perp) = \Gamma(L)$ for each $\v\in V$ and i) is proved. The last assertion follows from the equivalence of i) and iii).
\qed

\section{The Abstract Morimoto Theorem}\label{sec:morimoto}

\begin{mydef}
Let $E_1',E_2'\in \E(M)$ such that $\langle E_1', E_2'\rangle=0$ and such that there exist framings $V_i\subseteq \Gamma(E_i')$. Given $\Psi\in \SGF(E'_1\oplus E_2')$ we say that the triple $(V_1,V_2,\Psi)$ is {\it admissible} if there exists an isomorphism $\phi:\Gamma(E'_1\otimes\C)\to\Gamma(E'_2\otimes \C)$ of $C^\infty(M,\C)$-modules such that
\begin{enumerate}[i)]
\item $\phi(V_1\otimes\C)=V_2\otimes\C$;
\item $L_\Psi=\Gamma_\phi$;
\end{enumerate}
where $L_\Psi$ is the $\sqrt{-1}$-eigenbundle of $\Psi$ and $\Gamma_\phi=\{\e+\phi(\e)\,|\,\e\in E_1'\otimes \C\}$ is the graph of $\phi$. If this is the case, we say that $\phi$ is an {\it admissible isomorphism} for the admissible triple $(V_1,V_2,\Psi)$.
\end{mydef}

\begin{example}\label{ex:S^3admissible}
Consider the product manifold $M=M_1\times M_2$ in which each factor is a copy of $S^3$. For $i=1,2$ we pick  global frames $\{X_1^i,X_2^i,X_3^i\}$ (resp.\ $\{\alpha_1^i,\alpha_2^i,\alpha_3^i\}$) of $TM_i$ (resp.\ of $T^*M_i$ and functions $h^i\in C^\infty(M_i,\C)$ defining split structures $E_i,E_i'\in \E(M_i)$ and generalized $F$-structures $J_i\in \SGF(E_i)$, as in Example \ref{ex:S^3}. Furthermore, let $V_i$ be framings of $E_i'$ as in Example \ref{ex:S^3normal}. Fix $\tau=a+\sqrt{-1}b\in \C\setminus \R$ and let $\Psi\in \SGF(E_1'\oplus E_2')$ be such that  $\Psi(X_1)=aX_1^1+bX_1^2$ and $\Psi(\x_1^2)=b\x_1^1-a\x_1^2$. If $\lambda=b/(a+\sqrt{-1})$, then $\phi$ defined by $\phi(X_1^1)=\lambda X_1^2$ and $\phi(\x_1^2)=-\lambda \x_1^1$ is an admissible isomorphism for the admissible triple $(V_1,V_2,\Psi)$. If $h^1=h^2=0$, then $\Psi$ is the complex structure of modulus $\tau$ on the elliptic fibers of the Calabi-Eckmann fibration $S^3\times S^3\to S^2\times S^2$ described in \cite{CE}.
\end{example}

\begin{rem}\label{rem:admissible}
Let $E_1',E_2'\in \E(M)$ be mutually orthogonal with global framings $V_i\subseteq \Gamma(E_i')$. Let $\Psi\in \SGF(E_1'\oplus E_2')\cap{\rm O}(V_1\oplus V_2)$ be such that  $\pi_{E'_2}\circ\Psi_{|V_1}:V_1\rightarrow V_2$ is invertible. Here $\pi_{E'_2}$ denotes the orthogonal projection onto $E'_2$. Under these assumptions, $(V_1,V_2,\Psi)$ is admissible. To see this, write
\[
\Psi=\left[\begin{array}{cc}A&B\\C&D\end{array}\right]
\]
with blocks corresponding to the decomposition $E'_1\oplus E'_2$. Admissibility implies that the maps $B,C$ are invertible, and that
\[
L_{\Psi}=\{\e-\sqrt{-1}A\e-\sqrt{-1}C\e:\e\in E'_1\otimes \C\}=\Gamma_{\phi}\, ,
\]
where $\phi=-B^{-1}(A-\sqrt{-1}{\rm Id})$. Note that in this case the admissible isomorphism $\phi$ is unique. Moreover, after a choice of orthonormal bases on $V_1$ and $V_2$ is made, the morphism $\Psi$ is uniquely represented as a matrix
\[
\Psi_0=\left[\begin{array}{cc}A_0&B_0\\C_0&D_0\end{array}\right]\in \mathfrak{o}(2l,2l)\cap{\rm O}(2l,2l)\,,
\]
where the admissibility translates into the condition $B_0,C_0\in{\rm GL}(2l,\R)$.In particular, the matrix
\[
\Psi_0^{\rm can}=\left[\begin{array}{cc}0&{\rm Id}\\-{\rm Id}&0\end{array}\right]
\]
yields the admissible triple used in the original work of Morimoto \cite{M} and in some of its generalizations \cite{N}, \cite{GT1}, \cite{GT2}.
\end{rem}

\begin{prop} 
Let $E_1',E_2'\in \E_l(M)$ be mutually orthogonal with global framings $V_i\subseteq \Gamma(E_i')$. Let $\Sigma\subseteq \SGF(E_1'\oplus E_2')\cap{\rm O}(V_1\oplus V_2)$ be the subset of all $\Psi$ such that $\pi_{E'_2}\circ\Psi_{|V_1}:V_1\rightarrow V_2$ is invertible. Then $\Sigma$ is homeomorphic to ${\rm O}(l,l)$.
\end{prop}

\noindent\emph{Proof}. The group ${\rm O}(V_1)\times {\rm O}(V_2)$ acts transitively on $\Sigma$ by conjugation or, more precisely, by
\[
(R_1,R_2)\cdot \Psi:= R\Psi R^{-1}\,,
\]
where $R=R_1\oplus R_2:V_1\oplus V_2\rightarrow V_1\oplus V_2$. Given $\Psi_0\in\Sigma$, Effros' Open Mapping Theorem \cite{Eff} shows that the canonical bijection
\[
\frac{{\rm O}(V_1)\times {\rm O}(V_2)}{{\rm Stab}(\Psi_0)}\rightarrow \Sigma
\]
defined by $(R_1,R_2){\rm Stab}(\Psi_0)\mapsto (R_1,R_2)\cdot \Psi_0$ is a homeomorphism. On the other hand, the stabilizer ${\rm Stab}(\Psi_0)$ consists of the pairs of the form $(R_1,\phi_0 R_1\phi_0^{-1})$ (where $\phi_0$ is the admissible isomorphism of $\Psi_0$), and the 
projection ${\rm O}(V_1)\times {\rm O}(V_2)\rightarrow {\rm O}(V_2)$ descends to a homeomorphism
\[
\frac{{\rm O}(V_1)\times {\rm O}(V_2)}{{\rm Stab}(\Psi_0)}\rightarrow {\rm O}(V_2)\,.
\]
Combining these observations, we obtain the following chain of homeomorphisms
\[
\Sigma\simeq \frac{{\rm O}(V_1)\times {\rm O}(V_2)}{{\rm Stab}(\Psi_0)} \simeq {\rm O}(V_2)\simeq {\rm O}(l,l)\,.
\]
\qed

\begin{rem}\label{rem:families}
If $l=1$ then $O(1,1)$ is one dimensional and the construction of Remark \ref{rem:admissible} yields a one-parameter family of admissible triples. A particular instance is the $\tau$-dependent family of admissible triples on $S^3\times S^3$ described in Example \ref{ex:S^3admissible}.
\end{rem}

\begin{lem}\label{lem:admissible}
Let $E_1',E'_2\in \E(M)$ be such that $\langle E_1',E_2'\rangle=0$ and let $V_i\subseteq \Gamma(E_i')\cap \I(E_i')$ be framings of $E_i'$. Given $\Psi\in \SGF(E_1'\oplus E'_2)$ such that $(V_1,V_2,\Psi)$ is an admissible triple, then $\Psi\in \CRF(E_1'\oplus E'_2)$ if and only if $\phi([\v,\w])=[\phi(\v),\phi(\w)]$ for all $\v,\w\in V_1$.
\end{lem}

\noindent\emph{Proof:} By assumption, $V_i\subseteq \I(E_i')$ and $\phi(V_1)\subseteq \I(E_2')$. Therefore, $[\v,\phi(\w)]\subseteq (E_1'\cap E_2')\otimes\C=0$ for any $\v,\w\in V_1$. It follows that $\Psi$ is integrable if and only if 
\[
0=\langle [\v+\phi(\v),\w+\phi(\w)],\z+\phi(\z)\rangle= \langle[\v,\w],\z\rangle + \langle[\phi(\v),\phi(\w)],\phi(\z)\rangle
\]
for every $\v,\w,\z\in V_1$. The isotropy of $\Gamma_\phi$ implies 
\[
\langle[\v,\w],\z\rangle=-\langle\phi([\v,\w]),\phi(\z)\rangle\,,
\]
which concludes the proof.
\qed

\begin{example}\label{ex:admissible}
If $\Psi$ is as in Remark \ref{rem:admissible} with $A=D=0$, then the admissible isomorphism $\phi$ maps $V_1$ to $\sqrt{-1}V_2$. If this is the case, Lemma \ref{lem:admissible} shows that $\Psi \in \CRF(E_1'\oplus E'_2)$ if and only if the Dorfman bracket vanishes when restricted to $V_1$ and $V_2$.
\end{example}

\begin{mydef}\label{def:morimoto}
Let $E_1,E_2,E_1',E_2'\in \E(M)$ be mutually orthogonal split structures. For $i=1,2$, denote $E_i''=E_i\oplus E_i'$ and let $E''=E_1''\oplus E_2''$.
A \emph{Morimoto datum} on $M$ is given by $(J_1,J_2,V_1,V_2,\Psi)$, where $J_i\in \SGF(E_i)$, $V_i$ is a framing of $E_i'$ for $i=1,2$ and $\Psi\in\SGF(E'_1\oplus E'_2)$, satisfies the following conditions:
\begin{enumerate}[1)]
\item $V_{i}\subseteq \I(E_{1}'')\cap\I(E_{2}'')$ for $i=1,2$;
\item there exist local framings $W_i\subseteq \I(E_1'')\cap \I(E_2'')$ of $E_i$ for $i=1,2$;
\item $(V_1,V_2,\Psi)$ is an admissible triple.
\end{enumerate}
We say that a Morimoto datum is {\it adaptable} if the local framings $W_i$ as above satisfy $d\langle J_i(W_i),W_i\rangle \subseteq \Gamma(E_i'')$. If such a $W_i$ exists, we call it an {\it adapted local framing} of $E_i$.
\end{mydef}

\begin{lem}\label{lem:adapted}
Let $\mathcal M=(J_1,J_2,V_1,V_2,\Psi)$ be a Morimoto datum.
\begin{enumerate}[i)]
\item If $J_i\in\CRF(E_i)$, then $\mathcal M$ is adaptable;
\item If $\mathcal M$ is adaptable, then $[\Gamma(L_{J_1}),\Gamma(L_{J_2})]=0$.
\end{enumerate}
\end{lem}

\noindent\emph{Proof:}
Let $W_i$ be local framings of $E_i$ as in Definition \ref{def:morimoto}. Since $J_i\in \CRF(E_i)$, $[\w-\sqrt{-1}J_i(\w),\z-\sqrt{-1}J_i(\z)]$ is in $\Gamma(E_i\otimes \C)$ for each $\w,\z\in W_i$. Taking the imaginary part, $[J_i(\w),\z]+[\w,J_i(\z)]$ is in  $\Gamma(E_i)$. Since $2d\langle J_i(\w),\z\rangle = [J_i(\w),\z]+[\z,J_i(\w)]$ and $[W_i,J_i(W_i)]\subseteq \Gamma(E_i'')$, we conclude that $W_i$ is an adapted local framing and thus i) holds. Let $W_1$ and $W_2$ be respective adapted local framings of $E_1$ and $E_2$. Notice that $[W_1,W_2]\in E_1''\cap E_2''=0$ and $[W_1,J_2(W_2)]\subseteq \Gamma(E_2'')$. On the other hand, for each $\x\in W_1$ and $\y,\z\in W_2$
\[
0=\L_\x\langle J_2(\y),\z\rangle = \langle \L_\x(J_2(\y)), \z\rangle + \langle J_2(\y), \L_\x(\z)\rangle =  \langle \L_\x(J_2(\y)), \z\rangle
\]
which implies $[W_1,J_2(W_2)]=0$. Similarly, $[J_1(W_1),W_2]=0$ and therefore $[J_1(W_1),J_2(W_2)]\in \Gamma(E_1\cap E_2)=0$. In particular, for each $\w_1\in W_1$ and $\w_2\in W_2$,
\[
[\w_1-\sqrt{-1}J_1(\w_1),\w_2-\sqrt{-1}J_2(\w_2)]=0\,.
\]
This concludes the proof since each $L_{J_i}$ is locally generated by sections of the form $\w_i-\sqrt{-1}J_i(\w_i)$.
\qed

\begin{lem}\label{lem:morimoto}
Let $\mathcal M=(J_1,J_2,V_1,V_2,\Psi)$ be a Morimoto datum. Then $(J_1,V_1)$ and $(J_2,V_2)$ are both normal pairs if and only if
\begin{enumerate}[1)]
\item $\mathcal M$ is adaptable;
\item $J_1$ and $J_2$ both normalize $J=J_1\oplus J_2\oplus \Psi$.
\end{enumerate}
\end{lem}

\noindent\emph{Proof:} Let $\Gamma=\Gamma(L_J)$ and let $\Gamma_i=\Gamma(L_{J_i})$ for $i=1,2$. Since the normality of the pair $(J_i,V_i)$ implies $J_i\in \CRF(E_i)$, Lemma \ref{lem:adapted} allows us to assume that $\mathcal M$ is adaptable and thus $[\Gamma_1,\Gamma_2]=0$. Since $V_i\subseteq \I(E_i'')$, we see that $[\Gamma_i,\Gamma_\phi]=[\Gamma_i,V_i]$ which implies that $[\Gamma_i,\Gamma]\subseteq \Gamma$ if and only if
\[
[\Gamma_i,\Gamma_i\oplus V_i]\subseteq \Gamma(L_J\cap E_i'')=\Gamma_i
\]
if and only if $(J_i,V_i)$ is a normal pair for $i=1,2$. \qed

\begin{theorem}[Abstract Morimoto Theorem]\label{Morimototheorem}
Let $\mathcal M=(J_1,J_2,V_1,V_2,\Psi)$ be a Morimoto datum. Then $\mathcal M$ satisfies
\begin{enumerate}[i)]
\item $J=J_1\oplus J_2\oplus \Psi$ is a generalized \CRF-structure;
\item $\mathcal M$ is adaptable;
\end{enumerate}
if and only if $\mathcal M$ satisfies
\begin{enumerate}[i')]
\item $(J_1,V_1)$ and $(J_2,V_2)$ are normal pairs;
\item $\Psi$ is a generalized \CRF-structure.
\end{enumerate}
\end{theorem}

\noindent\emph{Proof:} If $(J_1,V_1)$ and $(J_2,V_2)$ are normal pairs and $\Psi\in\CRF(E_1'\oplus E_2')$, then $J\in \CRF(E'')$ and $\mathcal M$ is adaptable by Lemma \ref{lem:morimoto}. Conversely, if $J\in \CRF(E'')$ then in particular $J_i$ normalizes $J_1\oplus J_2\oplus \Psi$. If in addition $\mathcal M$ is adaptable, then Lemma \ref{lem:morimoto} implies that $(J_1,V_1)$ and $(J_2,V_2)$ are normal pairs. As a consequence of Lemma \ref{lem:normality}, $V_i\subseteq \I(E_i')$ for $i=1,2$. Therefore, the admissible triple $(V_1,V_2,\Psi)$ satisfies the assumptions of Lemma \ref{lem:admissible} and therefore $\Psi$ is a generalized \CRF-structure.\qed

\section{Flat Ehresmann connections}\label{sec:connection}

\noindent
In this section we consider a surjective submersion $\pi:N\rightarrow M$ equipped with a flat Ehresmann connection, i.e.\ an involutive subbundle $H\subseteq TN$ such that
\[
TN=H\oplus \ker(T\pi)\,.
\]
The connection induces a splitting
\[
\T N=\left(H\oplus {\rm Ann}(\ker(T\pi))\right)\oplus\left(\ker(T\pi)\oplus{\rm Ann}(H)\right)\,.
\]
We refer to the split structures
$H\oplus {\rm Ann}(\ker(T\pi))$ and
$\ker(T\pi)\oplus{\rm Ann}(H)$, respectively, as the \emph{horizontal} and \emph{vertical split structure} defined by the connection $H$.

\begin{rem}\label{rem:PBbracket}
There is a canonical orthogonal isomorphism between
\[
\pi^*\T M=\{(q,X_p+\alpha_p): X_p+\alpha_p\in \T_pM, p\in M, q\in N, \pi(q)=p\}
\]
and $H\oplus {\rm Ann}(\ker(T\pi))$ given by the map
\[
(q,X_p+\alpha_p)\mapsto \hat{X}_q+\alpha_p\circ T_q\pi\, ,
\]
where $\hat{X}_q\in H_q$ is uniquely defined by $(T_q\pi)(\hat{X}_q)=X_{\pi(q)}$. Under this identification, $\pi^*\x\in \Gamma(\pi^*\T M)$ is the \emph{horizontal lifting} of $\x\in \Gamma(\T M)$. In particular, the restriction of $\pi^*$ to $\Gamma(T^*M)$ coincides with the usual pull-back of forms.
\end{rem}

\begin{lem}\label{lem:PBbracket}
For all $\x,\y\in\Gamma(\T M)$, $[\pi^*\x,\pi^*\y]=\pi^*[\x,\y]$.
\end{lem}

\noindent
\emph{Proof}: If $\x,\y$ are both forms, then both commutators vanish. If $\x,\y$ are both vector fields, the identity is a consequence of flatness. By linearity of the Dorfman bracket, it remains to consider the case $\x=\alpha\in \Gamma(T^*M)$ and $\y=X\in \Gamma(TM)$. Since 
\[
2\langle[\pi^*\alpha,\pi^*X], Y\rangle=(d\pi^*\alpha)(Y,\pi^*X)=(d\alpha)(T{\pi}Y,X)\circ\pi=2\langle\pi^*[\alpha,X], Y\rangle
\]
for all $Y\in\Gamma(TN)$, this shows that $[\pi^*\alpha,\pi^*X]=\pi^*[\alpha,X]$. Together with
$d\langle\pi^*X,\pi^*\alpha\rangle=\pi^*d\langle X,\alpha\rangle$, this concludes the proof.\qed

\begin{lem}\label{lem:hornormalizesver}
If $\x\in\Gamma(\T M)$, then $\pi^*\x\in\I(\ker(T\pi)\oplus {\rm Ann}(H))$.
\end{lem}

\noindent
\emph{Proof:} If $\v\in\Gamma((\ker(T\pi)\oplus {\rm Ann}(H))$ and $\x,\y\in\Gamma(\T M)$, then
\[
\langle[\pi^*\x,\v],\pi^*\y\rangle=-\langle\v,[\pi^*\x,\pi^*\y]\rangle=0\,.
\]
\qed

\begin{prop}\label{prop:PBnormalizer}
Let $S$ be a (real or complex) subbundle of $\T M\otimes \C$ and let $\x\in \Gamma_U(E)$, for some open set $U\subseteq M$. Then, $\x\in \I_U(E)$ if and only if $\pi^*\x\in \I_{\pi^{-1}(U)}(\pi^*E)$.
\end{prop}

\noindent\emph{Proof}: Let $\x\in \I_U(E)$, and let $U'\subseteq U$ be any open set that trivializes $S$. Given a frame $\{\v_i\}$ of $E$ on $U'$, then for all $\w\in \Gamma_{\pi^{-1}(U)}(\pi^*E)$, we have $\w_{\pi^{-1}(U')}=\sum_if_i\pi^*\v_i$ for some smooth functions $f_i$ defined on $\pi^{-1}(U')$, so that by Lemma \ref{lem:PBbracket}
\[
[\pi^*\x,\w]_{\pi^{-1}(U')}=\sum_i[\pi^*\x,f_i\pi^*\v_i]_{\pi^{-1}(U')}\in \Gamma_{\pi^{-1}(U')}(\pi^*(E))\,.
\]
Here and below, $[-,-]_O$ denotes the restriction of the Dorfman bracket to an open set $O$. Since the open sets $\pi^{-1}(U')$ cover $\pi^{-1}(U)$, we obtain
$[\pi^*\x,\w]\in \Gamma_{\pi^{-1}(U)}(\pi^*(E))$ and thus $\pi^*\x\in \I_{\pi^{-1}(U)}(\pi^*(E))$. Conversely, suppose that $\pi^*\x\in \I_{\pi^{-1}(U)}(\pi^*(E))$ and let $U'\subseteq U$, $\{\v_i\}$ be as before. If $\z\in \Gamma_U(E)$, then
$[\pi^*\x,\pi^*\z]\in \Gamma_{\pi^{-1}(U)}(\pi^*E)$. On the other hand, 
$[\pi^*\x,\pi^*\z]_{\pi^{-1}(U')}=\sum_i g_i\pi^*\v_i$ and from Lemma \ref{lem:PBbracket} we obtain
\[
\pi^*([\x,\z]_{U'})=\sum g_i\pi^*\v_i\,.
\]
It follows that $g_i=\pi^*h_i$, where $h_i$ are smooth functions $U'$ and
\[
[\x,\z]_{U'}=\sum h_i\v_i\in \Gamma_{U'}(E).
\]
Therefore,  $[\x,\z]\in \Gamma_{U}(E)$ and the proof is complete.
\qed

\begin{cor}\label{cor:PBnormal}
Let $E,E'$ be orthogonal split structures on $M$, let $J\in\SGF(E)$ and let $V$ be a framing of $E'$. Then $(J,V)$ is a normal pair if and only if $(\pi^*J,\pi^*V)$ is a normal pair.
\end{cor}

\section{Morimoto products}\label{sec:externalmorimoto}

\noindent For the remainder of the section we fix a  product manifold $N=M_1\times M_2$. In this case, we have submersions  $\pi_i: N\rightarrow M_i$
given by the projections onto the two factors. As in Remark \ref{rem:PBbracket} we obtain flat connections $H_i:=\ker(T\pi_j)$ and canonical isomorphisms
\[
\pi_i^*(\T M_i)\cong H_i\oplus {\rm Ann}(H_j)
\]
for $i\neq j$. Let us fix orthogonal split structures $E_i,E_i'\in \E(M_i)$, framings $V_i$ of $E_i'$ and split generalized F-structures $J_i$ on $E_i$. We also define $E_i''=E_i\oplus E_i'$, $E''=E_1''\oplus E_2''$ as well as
\[
\tilde E_i:=\pi_i^* E_i\,;\quad \tilde E'_i:=\pi_i^*E'_i\,;\quad\tilde J_i:=\pi^*_i J_i\,;\quad \tilde V_i:=\pi^*_i V_i\,.
\]
Note that $\tilde E_1,\tilde E_2,\tilde E_1',\tilde E_2'$ are mutually orthogonal split structures on $N$, $\tilde V_i$ is a framing of $\tilde E_i'$ and $\tilde J_i\in\SGF(\tilde E_i)$.

\begin{mydef}\label{def:externalmorimoto}
Let $E_i,E_i',J_i,V_i$ as above. Then $(J_1,J_2,V_1,V_2,\Psi)$ is an \emph{external Morimoto datum on $N$} is given by if $\Psi\in\SGF(E'_1\boxplus E'_2)$ satisfies the following conditions:
\begin{enumerate}[1)]
\item $V_{i}\subseteq \I(E_{i}'')$ for $i=1,2$;
\item  there exist local framings $W_i\subseteq \I(E_i'')$ of $E_i$ for $i=1,2$;
\item $(\pi_1^*V_1,\pi_2^*V_2,\Psi)$ is an admissible triple.
\end{enumerate}
An external Morimoto datum is called \emph{adaptable} if the local framings $W_i$ of condition 2) additionally satisfy $d\langle W_i,J_iW_i\rangle\subseteq \Gamma( E_i'')$ for $i=1,2$.
\end{mydef}

\begin{rem}\label{rem:externaladaptable}
If $E'_1=E_1^\bot$ and $E'_2=E_2^\perp$, then conditions 1) and 2) in Definition \ref{def:externalmorimoto} are trivially satisfied. Moreover, in this case all Morimoto data are adaptable.
\end{rem}

\begin{lem}\label{lem:externalmorimoto}
If $(J_1,J_2,V_1,V_2,\Psi)$ is an (adaptable) external Morimoto datum, then $(\pi^*_1J_1,\pi^*_2J_2,\pi^*_1V_1,\pi^*_2V_2,\Psi)$ is an (adaptable) Morimoto datum.
\end{lem}

\noindent\emph{Proof:} Let $(J_1,J_2,V_1,V_2,\Psi)$ be an external Morimoto datum, and let $\v\in V_i$. By Proposition \ref{prop:PBnormalizer} and Lemma \ref{lem:hornormalizesver}, $V_i\subseteq \I(\pi_1*E_1'')\cap \I(\pi_2^*E_2'')$.
Similarly, if $W_i$ is a local framing of $E_i$ as in Definition \ref{def:externalmorimoto}, then $\pi_i^*W_i$ is a local framing of $\pi_i^*E_i$ such that $\pi_i^*W_i\subseteq\I(\pi_1^*E_1'')\cap \I(\pi_2^*E_2'')$. The adaptability of $W_i$ implies
$d\langle\pi_i^*J_i\pi_i^*W_i,\pi_i^*W_i\rangle=\pi_i^*d\langle J_iW_i,W_i\rangle\subseteq \Gamma(\pi^*_iE_i'')$, which concludes the proof.
\qed

\begin{mydef}
Let $(J_1,J_2,V_1,V_2,\Psi)$ be an external Morimoto datum for $M_1\times M_2$. We define the \emph{Morimoto product of $J_1$ and $J_2$ with respect to $\Psi$} to be 
\[
J_1\boxplus_{\Psi}J_2:=\pi_1^*J_1\oplus\pi_2^*J_2\oplus\Psi\in \SGF(E'')\,.
\]
\end{mydef}

\begin{theorem}\label{thm:externalmorimoto}
Let $(J_1,J_2,V_1,V_2,\Psi)$ be an adaptable external Morimoto datum for $M_1\times M_2$.
Then $J_1\boxplus_{\Psi}J_2\in \CRF(E'')$ if and only if
\begin{enumerate}[i)]
\item $(J_1,V_1)$ and $(J_2,V_2)$ are normal pairs;
\item $\Psi\in \CRF(E_1'\boxplus E_2')$.
\end{enumerate}
\end{theorem}

\noindent\emph{Proof:} The result is a direct consequence of the Abstract Morimoto Theorem, which can be applied because of Lemma \ref{lem:externalmorimoto} and Corollary \ref{cor:PBnormal}.
\qed

\begin{rem}
If $J_i$ are generalized almost complex structures, then $V_1=V_2=0$ and thus $\Psi=0$. In this case, Theorem \ref{thm:externalmorimoto} amounts to the assertion that  $J_1\boxplus J_2$ is integrable if and only if both $J_1$ and $J_2$ are integrable.
\end{rem}

\begin{cor}\label{cor:externalmorimoto}
Let $J_i$ be generalized almost contact structures, let $E_i'=E_i^{\bot}$ and let $(V_1,V_2,\Psi)$ be admissible. Then $J_1\boxplus_{\Psi} J_2$ is a generalized complex structure on $M$ if and only if $(J_1,V_1)$ and $(J_2,V_2)$ are normal pairs, i.e.\ the generalized almost contact triples associated with $(J_i,V_i)$ in Example \ref{ex:GACS} are normal.
\end{cor}

\noindent\emph{Proof}: It suffices to observe that since ${\rm dim}\ V_i=2$, the normality of $(J_i,V_i)$ implies that the Dorfman bracket vanishes identically on $V_i$. Therefore, the admissible isomorphism $\phi$ satisfies
\[
[\phi(\v),\phi(\w)]=0=\phi[\v,\w]
\]
for all $\v,\w\in V_i$.
\qed

\begin{rem}
In particular, the integrability of Morimoto products of generalized almost contact structures does not depend on the choice of admissible triple.
\end{rem}

\begin{example}\label{ex:S^3morimoto}
If $M_1=M_2=S^3$ and $\tau,J_i,V_i,\Psi$ are as in Example \ref{ex:S^3admissible}, then $(J_1,J_2,V_1,V_2,\Psi)$ is an adaptable external Morimoto datum on $M=M_1\times M_2$. According to Corollary
\ref{cor:externalmorimoto} and Example \ref{ex:S^3normal}, $J=J_1\boxplus_\Psi J_2$ is integrable if and only if $h^i\in \ker(\overline\partial^i)\cap \ker(Y_1^i)$ for $i=1,2$. If $h^1=h^2=0$, these conditions are trivially satisfied and $J$ coincides with the family (parametrized by $\tau$) of complex structures on $S^3\times S^3$ discovered in \cite{CE}. On the other hand if $J$ is integrable and $(h^1,h^2)\neq 0$, then $J$ is a generalized complex structure that preserves $TM$ but not $T^*M$. As observed in \cite{K}, this implies that turning on the parameters $h^i$ has the effect of deforming the complex structure of Calabi and Eckmann by means of a holomorphic Poisson bivector. Therefore, the Morimoto product of two of the normal generalized almost contact structures on $S^3$ described in \cite{AG} with respect to split generalized $F$-structures $\Psi$ introduced in Example \ref{ex:S^3admissible} is a (generically non-commutative) Calabi-Eckmann structure on $S^3\times S^3$.

\end{example}

\section{Products with Lie Groups}\label{sec:sekiya}

\noindent For the remainder of this section let us fix a finite-dimensional Lie group $G$ with identity $e$ and a manifold $M$. We denote by $\mathfrak{g}$ the Lie algebra of $G$ we fix a basis $\{b_i\}$ of $\mathfrak{g}\ltimes\mathfrak{g}^*=\T_eG$. We also consider the left-action of $G$ acts on $M\times G$ defined by $h(p,g):=(p,hg)$ for all $p\in M$ and $g,h\in G$.

\begin{theorem}\label{thm:sekiya}
The following sets are in canonical bijection
\begin{enumerate}[i)]
\item $G$-invariant generalized almost complex structures $\mathcal J$ on $M\times G$, such that $\pi_{\T M}{\mathcal J}_{|\T G}$ is fiberwise injective, with image of split signature;
\item quadruples $(E,J,\{\v_i\},\varphi)$, where $E\in \E(M)$, $J\in \SGF(E)$ , $\{\v_i\}$ is a global frame of $E^{\bot}$ and $\varphi:M\rightarrow \mathfrak{o}(\mathfrak{g}\ltimes\mathfrak{g}^*)$ is a smooth map such that, for all $p\in M$
$$
\langle(\varphi^2_p+{\rm Id}_{\mathfrak{g}\ltimes\mathfrak{g}^*})b_i,b_j\rangle= \langle\v_i(p),\v_j(p)\rangle\,.
$$
\end{enumerate}
\end{theorem}

\noindent \emph{Proof:}
Given $\mathcal J$, for all $(p,g)\in M\times G$ we have
$$
{\mathcal J}_{p,g}=\left[\begin{array}{cc}A_{p,g}&B_{p,g}\\C_{p,g}&D_{p,g}\end{array}\right]\,,
$$
with respect to the decomposition $\T_{(p,g)}(M\times G)=\T_p M\oplus \T_gG$.
If
$\v_i(p):=B_{p,e}(b_i)$, then $E':=\Span(\{\v_i\})\subseteq\T M$ is a split structure and so is  $E:=(E')^{\bot}$. Let $J$ be defined by $J_p:=\left.A_{p,e}\right|_{E_p}$ for each $p\in M$. Since $E_p=\ker(C_{p,e})$, $J\in \SGF(E)$ and $\varphi$ defined by $\varphi_p:=D_{p,e}$ for each $p\in M$ is the required map. Conversely, consider a quadruple $(E,\{\v_i\},J,\varphi)$. For each $p\in M$, define
\[
\Psi_{p,e}:E_p^{\bot}\oplus\T_eG\rightarrow E_p^{\bot}\oplus\T_eG
\]
such that
\[
\Psi_{p,e}=\left[\begin{array}{cc}-B_{p,e}\varphi_pB_{p,e}^{-1}&B_{p,e}\\-B_{p,e}^*&\varphi_p\end{array}\right]\,,
\]
where
$B_{p,e}:\T_eG\rightarrow\T_pM$ is the isomorphism defined by $B_{p,e}(b_i):=\v_i(p)$. This map extends uniquely to a $G$-invariant bundle endomorphism $\Psi\in \SGF(E^{\bot}\boxplus\T G)$.
Let $\pi_1,\pi_2$ be the projections of $M\times G$ onto the respective factors. If $V_1=\Span_{\R}(\{\v_i\})$ and $V_2$ is the space of left-invariant sections of $\T G$, then
$(\pi_1^*V_1,\pi_2^*V_2,\Psi)$ is an admissible triple which gives rise to the Morimoto product
\[
{\mathcal J}:=J\boxplus_\Psi 0\,.
\]
The assignments
${\mathcal J}\mapsto (E,\{\v_i\},J,\varphi) $ and $(E,\{\v_i\},J,\varphi)\mapsto  {\mathcal J}$ just described provide the required canonical bijections.
\qed

\begin{example}
If $G=\R$ the condition $\langle(\varphi^2_p+{\rm Id})b_i,b_j\rangle= \langle\v_i(p),\v_j(p)\rangle$ and, correspondingly, the condition that $\pi_{\T M}\mathcal J(\T\R)$ is of split signature are both automatically satisfied. Therefore, Theorem \ref{thm:sekiya} reduces to Sekiya's characterization \cite{S} of invariant generalized almost complex structure on $M\times \R$.
\end{example}

\begin{cor}
Let $\mathcal J$ be a
$G$-invariant generalized almost complex structure on $M\times G$, such that $B:=\pi_{\T M}{\mathcal J}_{|\T G}$ is fiberwise injective, with image of split signature. Let  $(E,J,\{\v_i\},\varphi)$ corresponding to $\mathcal J$ under the bijection of Theorem \ref{thm:sekiya}. Then, $\mathcal J$ is integrable if and only if $(J,\Span_{\R}(\{\v_i\}))$ is a normal pair and the map
$$
\phi=(B^*)^{-1}\circ(\varphi-\sqrt{-1}{\rm Id})
$$
satisfies $[\phi(\v),\phi(\w)]=\phi[\v,\w]$ for all $\v,\w\in \Gamma(\T G)$.
\end{cor}

\noindent\emph{Proof:} As in the proof of Theorem \ref{thm:sekiya}, $\mathcal J$ can be written as a Morimoto product of the form $J\boxplus_{\Psi} 0$. By Remark \ref{rem:externaladaptable}, the corresponding Morimoto datum is adaptable. Theorem \ref{thm:externalmorimoto} then guarantees that $\mathcal J$ is integrable if and only if $(J,\Span_{\R}(\{\v_i\}))$ is a normal pair and $\Psi$ is a split generalized \CRF-structure. The result then follows from Lemma \ref{lem:admissible}.
\qed

\begin{example}\label{ex:R^k}
Let $G=\R^k$, and assume that $\mathcal J$ satisfies the conditions of Theorem \ref{thm:sekiya}. Then Example \ref{ex:admissible} shows that $\mathcal J$ is integrable if and only if $(J,\Span_{\R}(\{\v_i\}))$ is a normal pair and $[\v_i,\v_j]=0$ for all $i,j$.
\end{example}

\begin{example}
Let $J$ be a split generalized $F$-structure defined by a classical $F$-structure on a manifold $M$ as in Example \ref{ex:fstructures}. Suppose that $J$, together with vectors $\{v_i\}\subseteq \Gamma(TM)$, endows $M$ with the structure of {\it $f$-manifold with complemented frame} in the sense of \cite{N}. Let $\Psi_{0}^{\rm can}$ be as in Remark \ref{rem:admissible} with respect to the basis consisting of the complemented frame $\{\v_i\}$ extended by the standard orthogonal basis of invariant sections of $\T\R^k$. By definition, $M$ is a {\it normal framed $f$-manifold} if the generalized almost complex structure $J\boxplus_{\Psi_0^{\rm can}} 0$ is integrable. By Example \ref{ex:R^k}, we see that $M$ is a normal framed $f$ manifold if and only if $(J,\Span_{\R}(\{\v_i\}))$ is a normal pair and $[\v_i,\v_j]=0$ for all $i,j$. In \cite{N}, Nakagawa proved the following generalization of Morimoto's Theorem: the Morimoto product $J_1\boxplus_{\Psi_0^{\rm can}} J_2$ of two framed $f$-manifolds $(M_1,J_1,\{\v_{1,i}\})$ and $(M_2,J_2,\{\v_{2,i}\})$ is integrable if and only if $M_1$ and $M_2$ are normal framed $f$-manifolds. Thanks to Example \ref{ex:admissible}, one may view Nakagawa's Theorem as a particular case of Theorem \ref{thm:externalmorimoto}.
\end{example}

\section{Flat principal bundles}\label{sec:flat}

\noindent
For the reminder of this section, let $\pi:N\rightarrow M$ be a principal bundle with fiber $G$ admitting a flat connection $H$. As customary in this context, we assume $H$ to be $G$-invariant, so that the vertical and horizontal split structures are $G$-invariant as well. If a basis $\{v_i\}$ of the Lie algebra of $G$ is fixed and $\tilde v_i\in \Gamma(TN)$ denotes the fundamental vector field generated by $v_i$, then $\ker(T\pi)$ is trivialized by the global frame $\{\tilde v_i\}$ while ${\rm Ann}(H)$ is trivialized by the dual global coframe $\{\tilde v_i^*\}$. In particular, the vertical split structure $\ker(T\pi)\oplus {\rm Ann}(H)$ is a trivial bundle. Moreover, the framing $V'=\Span_{\R}(\{\tilde v_i,\tilde v_j^*\}_{i,j})$ of the vertical split structure is involutive, i.e.\ it is closed under the Dorfman bracket.

\begin{lem}\label{lem:flat}
Let $E,E'$ be orthogonal split structures on $M$, let $J\in \SGF(E)$ and let $V$ be a framing of $E'$. If $E''=E\oplus E'$, then
\begin{enumerate}[i)]
\item $V'\subseteq \I(\pi^*E'')\cap\I(\ker(T\pi)\oplus {\rm Ann}(H))$;
\item $\pi^*V\subseteq \I(\pi^*E'')\cap\I(\ker(T\pi)\oplus {\rm Ann}(H))$ if and only if $V\subseteq \I(E'')$.
\end{enumerate}
\end{lem}

\noindent\emph{Proof:}
In order to prove the first statement, let $\v,\w\in V'$ and let $\u\in\Gamma(\T M)$. From the involutivity of $V'$, it follows that $\langle[\v,\w],\pi^*\u\rangle=0$ and thus $V'$ normalizes the vertical split structure. Similarly, if $\e\in \Gamma(E'')$ then
\[
\langle [\v,\pi^*\e],\w\rangle=-\langle \pi^*\e,[\v,\w]\rangle=0\,.
\]
This, together with
\[
\langle [\v,\pi^*\e],\pi^*\u\rangle=\langle \v ,[\pi^*\e,\pi^*\u]\rangle=0\,,
\]
shows that $[\v,\pi^*\e]=0$ and thus $V'$ normalizes $\pi^*E''$.
The second statement is a direct consequence of Lemma \ref{lem:hornormalizesver} and Proposition \ref{prop:PBnormalizer}. \qed

\begin{theorem}\label{thm:flat} 
In addition to the assumptions of Lemma \ref{lem:flat}, suppose that
\begin{enumerate}[i)]
\item $V\subseteq \I(E'')$;
\item $E'$ admits local framings such that $W\subseteq\I(E'')$ and $d\langle W,JW\rangle\subseteq \Gamma(E'')$;
\item $(\pi^*V,V', \Psi)$ is an admissible triple.
\end{enumerate}
Then $\pi^*J\oplus\Psi\in \CRF(\pi^*E''\oplus \ker(T\pi)\oplus {\rm Ann}(H))$ if and only if
$(J,V)$ is a normal pair and the admissible isomorphism $\phi:V'\rightarrow \pi^*V$ is a Lie algebra isomorphism.
\end{theorem}

\noindent\emph{Proof:} Our assumptions, together with Lemma \ref{lem:flat} imply that $(\pi^*J, 0, \pi^*V, V',\Psi)$ is an adaptable Morimoto datum. The result then follows combining the Abstract Morimoto Theorem, Corollary \ref{cor:PBnormal} and Lemma \ref{lem:admissible}.
\qed

\begin{rem} Note that any flat principal $G$-bundle on $M$ can be written in the form $N=(\widetilde{M}\times G)/\pi_1(M)$, where $\widetilde{M}$ is the universal cover of $M$ and $\pi_1(M)$ acts on $G$ by holonomy. This point of view suggests an alternative method to construct split generalized $F$-structures on $N$. Start from a structure on $M$, lift it to a $\pi_1(M)$-invariant structure on $\widetilde{M}$, take a Morimoto product with a $\pi_1(M)$-invariant structure on $G$ and descend the resulting structure to $N$. In the context of classical contact geometry, this (in the more general context of flat bundles) is described in \cite{BH}. This should be contrasted with Theorem \ref{thm:flat} in which $G$ is not endowed with split generalized $F$ structures and instead an admissible triple is used to extend the \SGF structure on $\pi^*\T M$ to a possible larger split structure.
\end{rem}

\section{Abstract Blair-Ludden-Yano Theorem}\label{sec:BLY}

\begin{mydef}
Let $E\in \E(M)$. We say that $(V,W)$ is a {\it split framing of $E$} if $V$ and $W$ are isotropic and $V\oplus W$ is a framing of $E$.
\end{mydef}

\begin{mydef}\label{def:contactdatum}
Let $E\in \E(M)$ and $E'\in \E_k(M)$ be mutually orthogonal. Given a maximal isotropic subbundle $L\subseteq E$ and a split framing $(V,W)$ for $E'$,
we say that $(L,V,W)$ is a {\it rank $k$ contact datum for $(E,E')$} if
\begin{enumerate}[1)]
\item $V\subseteq \I(E)$;
\item $\Gamma(L)\oplus W\subseteq \I(\Gamma(L)\oplus W)$;
\item $[V\oplus W,V\oplus W]=0$;
\item $\Gamma(L)=\L_W(\Gamma(L))$;
\item $\Gamma(E)=\Gamma(L)\oplus \L_V(\Gamma(L))$.
\end{enumerate}
\end{mydef}

\begin{rem}\label{rem:contact1}
If $E'=E^\perp$, then 3) implies $\langle [V\oplus W,\Gamma(E)],V\oplus W\rangle=0$. In turn, this shows that condition 1) is automatically satisfied and that condition 2) simplifies to $\Gamma(L)\subseteq \I(\Gamma(L)\oplus W)$.
\end{rem}

\begin{rem}\label{rem:contact2}
If $E'\in \E_1(M)$ a split framing $(V,W)$ is uniquely determined by $V$. This observation allows us to use the shorthand notation $(L,V)$ for a rank 1 contact datum $(L,V,W)$.
\end{rem}

\begin{lem}\label{lem:contact}
If $(L,V)$ is a rank 1 contact datum for $(E,E')$, then
\begin{enumerate}[i)]
\item $\L_V(L)\subseteq E$ is maximal isotropic;
\item $W\subseteq \I(E)$.
\end{enumerate}
\end{lem}

\noindent\emph{Proof:}
If $\e$ is a generator of $V$ and $\x,\y\in \Gamma(L)$, then
\[
\langle \L_\e\x,\L_\e\y \rangle = \langle \L_\x\e,\L_\y\e \rangle = \L_\y\langle \L_\x\e,\e \rangle - \langle \e,\L_\x\L_\y\e \rangle = - \langle \e,\L_\x\L_\y\e \rangle \,.
\]
Similarly,
\[
\langle \L_\e\x,\L_\e\y \rangle = - \langle \L_\y\L_\x\e,\e \rangle = \langle \L_\x\L_\y\e,\e \rangle-\langle \L_{[\y,\x]}\e,\e  \rangle = \langle \L_\x\L_\y\e,\e \rangle
\]
from which i) follows. To prove ii) observe that for each $\w\in W$
\[
\langle \L_\w\x,\e\rangle = \L_\w\langle \x,\e\rangle - \langle \x,\L_\w\e \rangle =0
\]
implies $W\subseteq \I(L)$. On the other hand,
\[
\L_\w(\L_\e\x)=\L_{[\w,\e]}\x-\L_\e(\L_\w\x)\in \L_V(L)
\]
shows that $W\subseteq \I(\L_V(L))$ which concludes the proof.
\qed

\begin{example}\label{ex:contact}
Consider a contact form $\eta$ on $M$ and a corresponding Reeb vector field $\xi$. If $E'=\Span(\xi,\eta)$ and $E=(E')^\perp$, then $(TM\cap E,{\rm span}(\eta))$ is a rank 1 contact datum.
\end{example}

\begin{example}\label{ex:S^3contact}
In the notation of Example \ref{ex:S^3}, let $L=\Span(X_2,X_3)$ and $V=\Span_\R(\x_1)$. Thanks to Remark \ref{rem:contact1}, $(L,V)$ is a rank 1 contact datum for $(E,E')$ if and only if
\[
0=[X_1,\x_1]={\rm Re}(Y_1(h))X_2+{\rm Im}(Y_1(h))X_3
\]
and hence if and only if $Y_1(h)=0$. If $h=0$ this is a particular case of Example \ref{ex:contact}. On the other hand, if $h\neq 0$ then $\x_1$ is no longer a 1-form and therefore the resulting contact datum is not defined by a classical contact structure.
\end{example}

\begin{mydef}
Let $(L,V)$ be a rank 1 contact datum for $(E,E')$ and let $J\in \SGF(E)$. We say that $(J,L,V)$ is a {\it normal contact datum for $(E,E')$} if $J(L\oplus \Span(W))\subseteq L\oplus \Span(W)$ and $(J,V\oplus W)$ is a normal pair.
\end{mydef}

\begin{rem}
If $E'=E^\perp$, then combining Lemma \ref{lem:contact}, Remark \ref{rem:contact1} and Lemma \ref{lem:normality} we see that $(J,V\oplus W)$ is a normal pair if and only if $J\in \CRF(E)$.
\end{rem}

\begin{example}\label{ex:normalcontact}
Let $\xi$ and $\eta$ be as in Example \ref{ex:contact} and let $\phi\in {\rm End}(TM)$ be such that $(\phi, \xi,\eta)$ is a classical almost complex structure. If $J$ denotes the split generalized $F$-structure induced by $\phi$ on $E$, then $(J,TM\cap E, {\rm span}_\R(\eta))$ is a normal contact datum if and only if $(\phi,\xi,\eta)$ is a normal almost contact structure.
\end{example}

\begin{example}\label{ex:S^3normalcontact}
Let $(L,V)$ be the rank 1 contact datum of Example \ref{ex:S^3contact} and let $J$ be as in Example \ref{ex:S^3}. Then $(J,L,V)$ is a normal contact datum if and only if $h\in \ker(\overline\partial)\cap\ker(Y_1)$.
\end{example}

\begin{mydef}\label{def:bicontactdatum}
Let $E,E_1',E_2'$ be mutually orthogonal split structures with $E_1'$ and $E_2'$ of rank 1. Given a maximal isotropic subbundle $L\subseteq E$ and split framings $(V_1,W_1)$ for $E_1$ and  $(V_2,W_2)$ for $E_2$, we say that $(L,V_1,V_2)$ is a {\it bicontact datum for $(E,E_1',E_2')$} if
\begin{enumerate}[1)]
\item $(L,V_1\oplus V_2, W_1\oplus W_2)$ is a rank 2 contact datum for $(E,E_1'\oplus E_2')$;
\item $L=L_1\oplus L_2$, where $L_1=\ker(\L_{V_2})\cap L$ and $L_2=\ker(\L_{V_1})\cap L$;
\item For $i=1,2$, $L_i$ admits a local framing $K_i\subseteq \I(E\oplus E_1'\oplus E_2')$ such that $[K_1,K_2]=0$.
\end{enumerate}
If $(L,V_1,V_2)$ is a bicontact datum for $(E,E_1',E_2')$, we denote $E_i=L_i\oplus \L_{V_i}(L_i)$ and $E_i''=E_i\oplus E_i'$ for $i=1,2$. We also write $E''=E_1''\oplus E_2''$.
\end{mydef}

\begin{rem}\label{rem:bicontact}
If $E''=\T M$, then the condition $K_i\subseteq \I(E\oplus E_1'\oplus E_2')$ is automatically satisfied. For instance, this happens if $L_1\oplus W_1$ and $L_2\oplus W_2$ define complementary transverse foliations of constant rank. In this case, the condition $[K_1,K_2]=0$ is also satisfied by choosing local framings $K_i$ that are pushed-forward from the corresponding leaves.
\end{rem}

\begin{example}\label{ex:bicontact}
Let $\eta_1,\eta_2\in \Gamma(T^*M)$ be such that $(\eta_1,\eta_2)$ is an ordinary bicontact structure i.e.\ such that there exist $k_1,k_2\in \mathbb Z$ with the property that $\eta_1\eta_2(d\eta_1)^{k_1}(d\eta_2)^{k_2}$ is a volume form. Let $\xi_1, \xi_2\in\Gamma(TM)\cap \ker(\L_{\eta_1})\cap \ker(\L_{\eta_2})$ be  such that  $\<\eta_i,\xi_j\>=\delta_{ij}$ and let $E_i'=\Span(\eta_i,\xi_i)$ for $i=1,2$. If $E=(E_1'\oplus E_2')^\perp$, then $(TM\cap E, \Span(\eta_1),\Span(\eta_2))$ is a bicontact datum.
\end{example}

\begin{rem}\label{rem:bicontactmorimoto}
Let. $M=M_1\times M_2$. If $(L_i,V_i)$ are rank 1 contact data on $M_i$ for $i=1,2$, then, arguing as in Section \ref{sec:externalmorimoto}, $(L_1\boxplus L_2,\pi_1^*V_1,\pi_2^* V_2)$ is a bicontact datum on $M$.
\end{rem}

\begin{example}\label{ex:S^3bicontact}
Let $M_1=M_2=S^3$. If $(L_i,V_i,W_i)$ are the rank 1 contact data described in Example \ref{ex:S^3contact}, then $(L_1\boxplus L_2,\pi_1^*V_1,\pi_2^* V_2)$ is a bicontact datum on $S^3\times S^3$. Unless $h^1$ and $h^2$ both vanish (in which case we recover the standard bicontact structure on $S^3\times S^3$ of \cite{BLY}), this bicontact datum does not define a classical bicontact structure.
\end{example}

\begin{lem}\label{lem:bicontact}
Let $(L,V_1,V_2)$ be a bicontact datum for $(E,E_1',E_2')$. Then
\begin{enumerate}[i)]
\item $E_1$ and $E_2$ are orthogonal split structures;
\item $(L_1,V_1)$ and $(L_2,V_2)$ are rank 1 contact data;
\item $V_i\oplus W_i\subseteq \I(E_1'')\cap \I(E_2'')$;
\item $L_i$ admits local framings $K_i \subseteq \I(E_1'')\cap \I(E_2'')$ for $i=1,2$.
\end{enumerate}
\end{lem}

\noindent\emph{Proof:} Choose generators $\e_1\in V_1$ and $\e_2\in V_2$.  By assumption $\L_{V_i}(L_i)$ has the same rank as $L_i$. Therefore, arguing as in the proof of Lemma \ref{lem:contact}, $\L_{V_i}(L_i)$ is maximal isotropic in $E_i$ and thus $E_1,E_2\in \E(M)$. Since for each $\x_i\in \Gamma(L_i)$
\[
\< \x_1, \L_{\e_2}\x_2\> = \L_{\e_2} \< \x_1,\x_2\> - \< \L_{\e_2}\x_1,\x_2\>=0
\]
 and similarly, using $[\e_1,\e_2]=0$,
\[
\< \L_{\e_1}\x_1,\L_{\e_2}\x_2\> = \L_{\e_2} \< \L_{\e_1}\x_1,\x_2\> - \< \L_{\e_2}\L_{\e_1} \x_1,\x_2\> = - \<\L_{\e_1}\L_{\e_2} \x_1, \x_2 \> =0
\]
we conclude that $E_1$ and $E_2$ are orthogonal. In order to show that $(L_i,V_i)$ is a contact datum, we only need to check condition 2) in Definition \ref{def:contactdatum} since the remaining conditions are consequences of the assumption that $(L,V_1\oplus V_2, W_1\oplus W_2)$ is a rank 2 contact datum. Observe that for each $\x,\y\in \Gamma(L_1)\oplus W_1$
\[
\L_{\e_2} [ \x,\y]=[\L_{\e_2}\x,\y]+[\x,\L_{\e_2}\y]=0
\]
and
\[
\<[\x, \y],\e_2 \> =  \L_\x \<\y,\e_2 \> - \<\y,\L_\x \e_2 \> =0\,.
\]
Since $\Gamma(L)\oplus W$ is closed under the Dorfman bracket and $\Gamma(L_1)\oplus W_1$ is the subspace of $\Gamma(L)\oplus W$ orthogonal to $\e_2$ and annihilated by $\L_{\e_2}$, we conclude that $\Gamma(L_1)\oplus W_1$ is also closed under the Dorfman bracket. Hence $(L_1,V_1)$ is a rank 1 contact datum and, by the same token, so is $(L_2,V_2)$. Since $\L_{e_1}(K_2\oplus \L_{\e_2}(K_2))=0$, iii) is proved if we show that $W_1$ normalizes $E_2''$. To see this, observe that $[W_1,K_2]\subseteq \Gamma(L)\cap \ker \L_{\e_1} = L_2$. This implies that $\L_{\e_2}[W_1,K_2]\subseteq E_2$ and thus $[W_1,\L_{\e_2}K_2]$ are sections of $E_2$. Since $[W_1,V_2\oplus W_2]=0$, this concludes the proof of iii). Let $K_1$ and $K_2$ be local framings as in Definition \ref{def:bicontactdatum}. Then $0=\L_{\e_2}[K_1,K_2]=[K_1,\L_{\e_2}K_2]$ and similarly $[\L_{\e_1}K_1,K_2]=0$ so that $0=\L_{\e_2}[\L_{\e_1}K_1,K_2]=[\L_{\e_1}K_1,\L_{\e_2}K_2]$. We conclude that $K_i\oplus \L_{\e_i}(K_i)$ are mutually commuting local framings of $E_1$ and $E_2$. \qed

\begin{mydef}\label{def:Hermitian}
Let $(L,V_1,V_2)$ be a bicontact datum for $(E,E_1',E_2')$ and let $J\in \CRF(E'')$. We say that $(J,L,V_1,V_2)$ is a {Hermitian bicontact datum for $(E,E_1',E_2')$} if
\begin{enumerate}[1)]
\item $V_1\oplus W_1\subseteq \I(J)$;
\item $J(V_1)=V_2$ and $J(W_1)=W_2$;
\item $J(\Gamma(L)\oplus \Span(W_1\oplus W_2))\subseteq \Gamma(L)\oplus \Span(W_1\oplus W_2)$.
\end{enumerate}
\end{mydef}

\begin{example}\label{ex:Hermitian}
Let $(\eta_1,\eta_2)$ be a bicontact structure on $M$ and let $E$, $E_1'$, $E_2'$ be as in Example \ref{ex:bicontact}. If in addition $M$ is endowed with a Hermitian structure $(\mathcal J,g)$, then $M$ is said \cite{BLY} to be a {\it Hermitian bicontact manifold} provided that there exist $\xi_1,\xi_2\in \Gamma(TM)$ infinitesimal automorphisms of the Hermitian structure such that $\mathcal J(\xi_1)=\xi_2$ provided that $\eta_i$ is dual to $\xi_i$ with respect to the metric $g$. As proved in \cite{BLY}, these assumptions imply that $(TM\cap E, \Span(\eta_1),\Span(\eta_2))$ is a bicontact datum. If $J=\mathcal J\oplus (-\mathcal J^*)$ then all conditions in Definition \ref{def:Hermitian} are met, except possibly for $V_1\subseteq \I(J)$ which is equivalent to the requirement that $d\eta_1$ is of bidegree $(1,1)$ with respect to $\mathcal J$.
\end{example}

\begin{example}\label{ex:S^3hermitianbicontact}
Let $(L_1\boxplus L_2,\pi_1^*V_1,\pi_2^*V_2)$ be the bicontact datum on $S^3\times S^3$ introduced in Example \ref{ex:S^3bicontact} and let $J_i\in \SGF(M_i)$ be as in Example \ref{ex:S^3}. If $\Psi$ is as in Example \ref{ex:S^3morimoto}, then $(J_1\boxplus_\Psi J_2,L_1\boxplus L_2,\pi_1^*V_1,\pi_2^*V_2)$ is a Hermitian bicontact datum for $(E_1\boxplus E_2,\pi_1^*E_1',\pi_2^*E_2')$.
\end{example}

\begin{lem}\label{lem:Hermitian}
If $(J,L,V_1,V_2)$ is a Hermitian bicontact datum for $(E,E_1',E_2')$, then
\begin{enumerate}[i)]
\item $J(E_1)\subseteq E_1$, $J(E_2)\subseteq E_2$ and $J(E_1'\oplus E_2')\subseteq E_1'\oplus E_2'$;
\item If $J_1$ (resp.\ $J_2$) is the restriction of $J$ to $E_1$ (resp.\ $E_2)$ and $\Psi$ denotes the restriction of $J$ to $E_1'\oplus E_2'$, then $(J_1,J_2,V_1\oplus W_1,V_2\oplus W_2,\Psi)$ is a Morimoto datum.
\end{enumerate}
\end{lem}

\noindent\emph{Proof:} Let $\e_1\in V_1$, $\e_2\in V_2$ be generators. Since $\e_1\in \I(J)$ and $\e_2=J(\e_1)$, then $\e_2\in \I(J)$ by Lemma \ref{lem:Ju}.  In particular, $J(\ker (\L_{\e_2}))\subseteq \ker (\L_{\e_2})$. Moreover, $J(V_1)=V_2$ together with the orthogonality of $J$ imply that $\x\in \Gamma(E'')$ is orthogonal to both $V_1$ and $V_2$ if and only if $J(\x)$ is. Since by assumption $J(L_1)\subseteq L\oplus \Span(W_1\oplus W_2)$ and $L_1$ is the subbundle of $L\oplus \Span(W_1\oplus W_2)$ orthogonal to both $V_1\oplus V_2$ and annihilated by $\L_{\e_2}$, we conclude that $J(L_1)\subseteq L_1$. Since $J$ commutes with $\L_{\e_1}$, this implies that $J(\L_{\e_1}(L_1)\subseteq \L_{\e_1}(L_1)$ and thus $J(E_1)\subseteq E_1$. Similarly, $J(E_2)\subseteq E_2$. From Lemma \ref{lem:bicontact} we see that $V_i\oplus W_i\subseteq \I(E_1'')\cap\I(E_2'')$ and that $L_i$ admits local framings $K_i \subseteq \I(E_1'')\cap \I(E_2'')$. This proves the lemma since $(V_1\oplus W_1,V_2\oplus W_2,\Psi)$ is by construction an admissible triple. \qed

\begin{mydef}
We refer to $(J_1,J_2,V_1\oplus W_1,V_2\oplus W_2,\Psi)$ as in Lemma \ref{lem:Hermitian} as {\it the Morimoto datum corresponding to the Hermitian bicontact datum $(J,L,V_1,V_2)$}.
\end{mydef}

\begin{mydef}
A Hermitian bicontact datum $(J,L,V_1,V_2)$ is adaptable if for $i=1,2$, $L_i$ admits a local framing $K_i$ such that
\begin{enumerate}[1)]
\item $K_1,K_2\in \I(E'')$;
\item $[K_1,K_2]=0$;
\item $d\<J\x_i,\L_{V_i}(\y_i)\>\in \Gamma(E_i'')$ for any $\x_i,\y_i\in K_i$.
\end{enumerate}
If $K_i$ satisfies the above conditions, we say that $K_i$ is {\it an adapted local framing of $L_i$}.
\end{mydef}

\begin{example}\label{ex:Hermitianadaptable}
ofLet $M$ be a Hermitian bicontact manifold as in Example \ref{ex:Hermitian}. Then $L_1\oplus \Span(\xi_1)$ and $L_2\oplus \Span(\xi_2)$ define transverse foliations which by Lemma \ref{lem:Hermitian} are preserved by $\mathcal J$. Therefore, $J$ induces almost complex structures $\phi_i$ on the leaves $S_i$ of $L_i\oplus \Span(\xi_i)$. If for $i=1,2$ we let $K_i$ be the push forward under the inclusion map of a local framing of $TS_i$, then $K_i$ is an adapted local framing  $L_i$. Therefore, if $d\eta_1$ is of bidegree $(1,1)$ then the Hermitian bicontact datum constructed in Example \ref{ex:Hermitian} is adaptable.
\end{example}

\begin{rem}
Generalizing Example \ref{ex:S^3hermitianbicontact}, let $(L_1\boxplus L_2,\pi_1^*V_1,\pi_2^*V_2)$ be as in Remark \ref{rem:bicontact} and let $(J_i,L_i,V_i)$ be normal contact data for $i=1,2$. Given $\Psi\in \CRF(E_1'\boxplus E_2')$ such that $\Psi(V_1)\subseteq V_2$ and $\Psi(W_1)\subseteq \Psi(W_2)$, then $(J_1\boxplus_\Psi J_2,L_1\boxplus L_2, \pi_1^*V_1,\pi_2^*V_2)$ is an adaptable Hermitian bicontact datum.
\end{rem}

\begin{theorem}[Abstract Blair-Ludden-Yano Theorem]
Let $(J,L,V_1,V_2)$ be an adaptable Hermitian bicontact datum and let $(J_1,J_2,V_1\oplus W_1,V_2\oplus W_2,\Psi)$ be the corresponding Morimoto datum. Then $(J_1,L_1,V_1)$ and $(J_2,L_2,V_2)$ are normal contact data.
\end{theorem}

\noindent\emph{Proof:} Let $K_1$ and $K_2$ be adapted local framings of $L_1$ and $L_2$, respectively. Since $L_i$ is maximal isotropic, combining Lemma \ref{lem:bicontact} with Lemma \ref{lem:contact}, we see that $\L_{V_i}(L_i)$ is also maximal isotropic. Therefore, $K_i\oplus \L_{V_i}(K_i)$ is an adapted local framing of $E_i''$. Since by assumption $J\in \CRF(E'')$, then the Abstract Morimoto Theorem implies that $(J_1,V_1\oplus W_1)$ and $(J_2,V_2\oplus W_2)$ are normal pairs. As shown in the proof of Lemma \ref{lem:Hermitian}, $J_i$ preserves $L_i$ and thus $L_i\oplus \Span(W_i)$. Therefore, $(J_i,L_i,V_i)$ is a normal contact datum for $(E_i,E_i')$.

\begin{cor}[\cite{BLY}]
Let $M$ be a Hermitian bicontact manifold with $d\eta_1$ of bidegree $(1,1)$. Then $M$ is locally the product of two normal contact manifolds.
\end{cor}

\noindent\emph{Proof:} Example \ref{ex:Hermitianadaptable} shows that the Hermitian bicontact datum of the Hermitian bicontact manifold $M$ is adaptable and thus $(J_1,L_1,V_1)$  and $(J_2,L_2,V_2)$ are normal contact data by the Abstract Blair-Ludden-Yano Theorem. By Corollary \ref{cor:PBnormal}, $(J_i,L_i,V_i)$ induce normal contact data on the leaves $S_i$ of $L_i\oplus \Span(\xi_i)$. As observed in Example \ref{ex:normalcontact}, this implies that each leaf inherits the structure of normal contact manifold.
\qed

\end{document}